\documentclass[12pt]{report} 
\usepackage{latexsym}
\usepackage{amsfonts}
\usepackage[dvips]{graphicx}
\usepackage{amsthm}
\usepackage{amsmath}
\usepackage{amssymb}
\usepackage[arrow,matrix,curve,cmtip,ps]{xy}

\usepackage{verbatim}

\newcommand{\vsk}{\vskip 12pt}

\newcommand{\subj}[1]{\begin{center} \bf{ #1 } \end{center}}

\begin{document}


\title{Letters relating to a Theorem of Mr. Euler, of the Royal 
Academy of Sciences at Berlin, and F.R.S. for correcting the Aberrations 
in the Object-Glasses of refracting Telescopes. 
\footnote{Delivered to the Royal Society of London on April 9th, 1752, November 23rd, 1752, 
and July 8th, 1753. Originally published in the {\it Philosophical Transactions of the Royal Society of London}, volume 48, 1754. 
A copy of the original text is available electronically at the Euler Archive, 
http://www.eulerarchive.org/. These letters comprise E210 in the Enestr\"om Index.}
}
\author{James Short, John Dolland, Leonhard Euler 
\footnote{Portions authored by Euler translated from the French by Erik R. Tou, Dartmouth College.}
}
\date{\today}
\maketitle


\bigskip\bigskip

\subj{I. A Letter from Mr. James Short, F.R.S. to Peter Daval, Esq; F.R.S.}

\subj{(Read April 9, 1752.)}

[Surrey-street, April 9, 1752.]

\vsk

Dear Sir,

\vsk

There is published, in the Memoirs of the Royal Academy at Berlin, for the 
year 1747, a theorem by Mr. Euler, in which he shews a method of making 
object-glasses of telescopes, in such a manner, as not to be affected by 
the aberrations arising from the different refrangibility of the rays of 
light; these object-glasses consisting of two {\it meniscus} lens's, with 
water between them.

\vsk

Mr. John Dollond, who is an excellent analyst and optician, has examined the 
said theorem, and has discovered a mistake in it, which arises by assuming 
an hypothesis contrary to the established principles of optics; and, in 
consequence of this, Mr. Dollond has sent me the inclosed letter, which 
contains the discovery of the said mistake, and a demonstration of it.

\vsk

In order to act in the most candid manner with Mr. Euler, I have proposed 
to Mr. Dollond to write to him, shewing him the mistake, and desiring to 
know his reasons for that hypothesis; and therefore I desire, that this 
letter of Mr. Dollond's to me may be kept amongst the Soceity's papers, 
till Mr. Euler has had a sufficient time to answer Mr. Dollond's letter to 
him. I am,

\bigskip

{\hskip 250pt SIR,}

\vsk

{\hskip 250pt Your most humble servant,}

\vsk

{\hskip 250pt James Short.}

\break


\subj{II. A letter from Mr. John Dollond to James Short, A.M.F.R.S. 
concerning a Mistake in M. Euler's Theorem for correcting the Aberrations 
in the Object-Glasses of refracting Telescopes.}

\subj{(Read Nov. 23, 1752.)}

\vsk

[London, March 11, 1752.]

\vsk

SIR,

\vsk

The famous experiments of the prism, first tried by Sir Isaac Newton, sufficiently 
convinced that great man, that the perfection of telescopes was impeded by the 
different refrangibility of the rays of light, and not by the spherical figure 
of the glasses, as the common notion had been till that time; which put the 
philosopher upon grinding concave metals, in order to come at that by reflexion, 
which he despair'd of obtaining by refraction. For, that he was satisfied of 
the impossibility of correcting that aberration by a multiplicity of refractions, 
appears by his own words, in his treatise of Light and Colours, Book I. Part 2. 
Prop. 3.

\vsk

``I found moreover, that when light goes out of air through several contiguous 
mediums, as through water and glass, as often as by contrary refractions it is 
so corrected, that it emergeth in lines parallel to those in which it was incident, 
continues ever after to be white. But if the emergent rays be inclined to the 
incident, the whiteness of the emerging light will by degrees, in passing on 
from the place of emergence, become tinged in its edges with colours."

\vsk

It is therefore, Sir, somewhat strange, that any body now-a-days should attempt 
to do that, which so long ago has demonstrated impossible. But, as so great 
a mathematician as Mr. Euler has lately published a theorem\footnote{Vide Memoires 
of the Royal Academy of Berlin for the Year 1747.} for making object-glasses, that 
should be free from the aberration arising from the different refrangibility of light, 
the subject deserves a particular consideration. I have therefore carefully examined 
every step of his algebraic reasoning, which I have found strictly true in every part. 
But a certain hypothesis in page 285 appears to be destitute of support either from 
reason or experiiment, though it be there laid down as the foundation of the whole 
fabrick. This gentleman puts $m:1$ for the ratio of refraction out of air into glass of 
the mean refrangible rays, and $M:1$ for that of the least refrangible. Also for 
the ratio of refraction out of air into water of the mean refrangible rays he 
puts $n:1$, and for the least refrangible $N:1$. As to the numbers, he makes 
$m=\frac{31}{20}$, $M=\frac{77}{50}$, and $n=\frac{4}{3}$; which so far 
answer well enough to experiments. But the difficulty consists in finding the 
value of $N$ in a true proportion to the rest.

\vsk

Here the author introduces the supposition above-mention'd; which is, that $m$ 
is the same power of $M$, as $n$ is of $N$; and therefore puts $n=m^a$, and 
$N=M^a$. Whereas, by all the experiments that have hitherto been made, the 
proportion will come out thus, $m-1:n-1::m-M:n-N$.

\vsk

The leters fixed upon by Mr. Euler to represent the radii of the four refracting 
surfaces of his compound object-glass, are $f$ $g$ $h$ and $k$, and the distance 
of the object he expresses by $a$; then will the focal distance be 
$=\frac{1}{n(\frac{1}{g} - \frac{1}{h}) + m(\frac{1}{f} - \frac{1}{g} + \frac{1}{h} 
- \frac{1}{k}) - \frac{1}{a} - \frac{1}{f} + \frac{1}{k}}$. Now, says he, it is 
evident, that the different refrangibility of the rays would make no alteration, 
either in the place of the image, or in its magnitude, if it were possible to 
determine the radii of the four surfaces, so as to have $n(\frac{1}{g}-\frac{1}{h}) 
+ m(\frac{1}{f} - \frac{1}{g} + \frac{1}{k} - \frac{1}{h}) = N(\frac{1}{g}-\frac{1}{h}) 
+ M(\frac{1}{f} - \frac{1}{g} + \frac{1}{k} - \frac{1}{h})$. And this, Sir, I shall 
readily grant. But when the surfaces are thus proportioned, the sum of the refractions 
will be $=0$; that is to say, the emergent rays will be parallel to the incident. For, 
if $n(\frac{1}{g}-\frac{1}{h}) + m(\frac{1}{f} - \frac{1}{g} + \frac{1}{h} - \frac{1}{k}) 
= N(\frac{1}{g}-\frac{1}{h}) + M(\frac{1}{f} - \frac{1}{g} + \frac{1}{h} - \frac{1}{k})$, 
then $n - N(\frac{1}{g} - \frac{1}{h}) + m - M(\frac{1}{f} - \frac{1}{g} + \frac{1}{h} 
- \frac{1}{k}) = 0$. Also if $n-N:m-M::n-1:m-1$, then $n-1(\frac{1}{g} - \frac{1}{h}) 
+ m-1(\frac{1}{f} - \frac{1}{g} + \frac{1}{h} - \frac{1}{k})=0$; or otherwise $n(\frac{1}{g} 
- \frac{1}{h}) + m(\frac{1}{f} - \frac{1}{g} + \frac{1}{h} - \frac{1}{k}) - \frac{1}{f} 
+ \frac{1}{k} = 0$; which reduces the denominator of the fraction expressing the focal 
distance to $\frac{1}{a}$. Whence the focal distance will be $=a$; or, in other words, 
the image will be the object itself. And as, in this case, there will be no refraction, 
it will be easy to conceive how there should be no aberration.

\vsk

And now, Sir, I think I have demonstrated, that Mr. Euler's theorem is intirely founded 
upon a new law of refraction of his own; but that, according to the laws discover'd 
by experiment, the aberration arising from the different refrangibility of light at the 
object-glass cannot be corrected by any number of refractions whatsoever. I am,

\bigskip

{\hskip 250pt SIR,}

\vsk

{\hskip 250pt Your most obedient humble servant,}

\vsk

{\hskip 250pt John Dollond.}

\break


\subj{III. Mr. Euler's Letter to Mr. James Short, F.R.S.}

\subj{(Read July 8, 1753.)}

\vsk

[Berlin, 19 June, 1752.]

\vsk

Sir,

\vsk

You have done me a great service, in having disposed Mr. Dollond to suspend 
the proposal of his objections against my objective lenses until I have had an 
opportunity to respond, and to you I am infinitely obliged. I take, therefore, the 
liberty of addressing to you my response to him, requesting that you, after having 
deigned to your examination, send it to him at your pleasure; and in case you judge 
this material worthy of the attention of the Royal Society, I ask that you communicate 
the details of the proofs in my theory, which I have given in this letter. In any case, 
I hope that Mr. Dollond will be satisfied, since I agree with him on the little 
success that one could achieve with my lenses, by handling them according to the 
ordinary manner.

\vsk

\indent I am honored to be, with the most perfect esteem,

\bigskip

{\hskip 250pt Sir,}

\vsk

{\hskip 250pt Your most humble, and}

\vsk

{\hskip 250pt most obedient servant,}

\vsk

{\hskip 250pt L. Euler.}
 
\break


\subj{IV. To Sir Mr. Dollond}

\subj{(Read July 8, 1753.)}

\vsk

[Berlin, June 15, 1752.]

\vsk

Sir, 

\vsk

Being very appreciative of the honor you have shown me on the subject of the objective lenses which I had proposed, I must also tell you quite honestly that I have met with the most sizable obstacles in the execution of your design, consisting of four faces, which must conform exactly to the proportions that I had found. However, having done these experiments many times, which seem to have been more successful, we have found that the interval between the two lenses of the red and violet rays is very small, in which case there could not be a single lens of the same focal distance. Nevertheless, I must confess that such a lens, even if it were perfectly executed on my principles, would have other defects which would make matters worse than even ordinary lenses; it is that such a lens would not admit a very small opening due to the large curves which one must give to the interior faces, so that when one gives an ordinary opening, the image would deflect most confusedly.

\vsk

Thus since you have taken pains, Sir, to experiment on such lenses, and having done said experiments,\footnote{Mr. Dollond, in his letter to Mr. Euler, here referred to, does not say that he had made any trials himself, but only he had understood that such had been made by others, without success.} I pray that you discern the flaws, which can give rise to the different refrangibility of light rays, that arise from too large an opening; for this effect you will have only to leave a very small opening. However, if my theory were sound, as to which I will soon have the honor of speaking, it would be advisable to remedy this flaw: it would be necessary to give up on the spherical figure which one usually gives to the faces of the lenses, and endeavor to give another figure, and I have remarked that the figure of a parabola is advantageous, for it admits a very considerable opening. Our scientist Mr. Lieberkuhn has worked with lenses for which the curvature of the faces decreases from the middle to the edges, and he has found very great advantages. For these reasons I believe that my theory does not suffer anything on this front.

\vsk

For the theory, I agree with you, sir, that supposing the ratio of refraction of one medium in another unspecified medium for the middle light rays is $m$ to $1$, and for the red light rays is $M$ to $1$, the ratio of $m-M$ to $m-1$ will always be nearly constant, and that this conforms to nature, as the great Newton has remarked. This ratio does not differ more than almost imperceptibly from my theory: because since I maintain that $M=m^{\alpha}$, and that $m$ differs usually very little from $1$, namely $m=1+\omega$; and since $M=m^{\alpha}=1+\alpha \log(m)$ arbitrarily close, and $\log(1+\omega)=\log(m)=\omega$, also arbitrarily close, I will get $m-M=1+\omega-1-\alpha\omega=(1-\alpha)\omega$, and $m-1=\omega$, therefore the ratio $\frac{m-M}{m-1}$ will be $=1-\alpha$, very nearly constant. Further, I conclude that the experiments from which the great Newton had drawn his ratio could not contradict my theory.

\vsk

Secondly, I agree also that if the ratio $\frac{m-M}{m-1}=Const.$ were shown rigorously, it would not suffice to remedy the defect which results from the different refrangibility of the light rays, regardless of how one places the various transparent media, and that the interval between the various foci would always hold a constant ratio with the entire focal distance of the lens. But it is precisely this consideration, which to me furnishes the strongest argument: the eye seems to me such a perfect optical machine, but does not sense in any way the different refrangibility of the light rays. However small the focal distance may be, the sensitivity is so great that the various foci, if there be more than one, would not fail to considerably trouble the vision. However, it is quite certain that a healthy eye does not sense the effect of the different refrangibility.

\vsk

The marvelous structure of the eye, and the various humours of which it is composed, infinitely confirms me in this sentiment. For if it acted only to produce a representation on the back of the eye, one humour alone would have been sufficient; and the Creator would surely not have employed more. Further, I conclude that it is possible to eliminate the effect of the different refrangibility of the light rays by an appropriate arrangement of several trasparent media; thus since this would not be possible if the formula $\frac{m-M}{m-1}=Const.$ were shown rigorously, I draw the conclusion that it does not conform to nature.

\vsk

But here is a proof directing my thesis: I conceive of various transparent media, $A$, $B$, $C$, $D$, $E$, $etc$. which differ equally among themselves by ratio of their optical density, so that the ratio of refraction of each in the following will be the same. That is to say therefore in the passage from the first into the second the ratio of refraction for the red light rays $=r:1$, and for the violets $=v:1$; which will be the same in the passage from the second into the third, the third into the fourth, the fourth into the fifth, and so on. Moreover, it is clear that in the passage from the first into the third the ratio will be $=r^2:1$ for the red rays, and $=v^2:1$ for the violets; similarly, in the passage from the first into the fourth, the ratios will be $r^3:1$ and $v^3:1$.

\vsk

Therefore if in the passage through an unspecified medium the ratio of refraction of the red rays is $=r^n:1$, then the violet rays will be $v^n:1$; all this perfectly conforms to the principles of the great Newton. We suppose $r^n=R$, and $v^n=V$, such that $R:1$, and $V:1$ express the ratios of refraction for the red and violet rays in an arbitrary passage: and having $n\cdot \log(r)=\log(R)$ and $n\cdot \log(v)=\log(V)$ we will have $\log(R):\log(r)=\log(V):\log(v)$, whence $\frac{\log(R)}{\log(V)}=\frac{\log(r)}{\log(v)}$. Or put $v=r^{\alpha}$, and since $\log(v)=\alpha\cdot \log(r)$, one has $\frac{\log(R)}{\log(V)}=\frac{1}{\alpha}$, whence $\log(V)=\alpha\cdot \log(R)$, and then $V=R^{\alpha}$. 

\vsk

We have seen therefore the foundation of the principle which I employed in my article, which to me appears again unwavering; however I subject myself to the decision of the illustrious Royal Society, and to your judgement in particular, having the honor to be with the most perfect esteem, Sir,

\bigskip

{\hskip 250pt Your most humble}

\vsk

{\hskip 250pt and most obedient servant,}

\vsk

{\hskip 250pt L. Euler.}

\end{document}